\documentclass[12pt]{article}
\usepackage{amssymb,pstricks,graphicx,times}

\renewcommand{\P}{\mbox{$\mathbb{P}$}}

\def\calG{{\cal G}}
\let\hat\widehat

\begin{document}
\baselineskip 24pt

\title{Improving power in\\ genome-wide association studies:\\ 
weights tip the scale}

\vskip 1.5in
\author{
Kathryn Roeder,  B Devlin, and Larry Wasserman \\
\\ \\ \\ \\
  Department of Statistics (KR, LW)\\
  Carnegie Mellon University\\
  Pittsburgh,PA \\
  Department of Psychiatry (BD)\\
  University of Pittsburgh \\
  Pittsburgh,PA \\
    }
\date{}
\maketitle

\vskip .5in

\noindent Running Title: Weights in genome-wide association studies\\

\vskip .5in

\noindent Address for correspondence and reprints:

\noindent Kathryn Roeder, Department of Statistics, Carnegie Mellon
University, 5000 Forbes Avenue, Pittsburgh, PA 15213. E-mail:
roeder@stat.cmu.edu \\

\centerline{\bf Abstract}

Genome-wide association analysis has generated much discussion about
how to preserve power to detect signals despite the detrimental effect
of multiple testing on power.  We develop a weighted multiple testing
procedure that facilitates the input of prior information in the form
of groupings of tests.  For each group a weight is estimated from the
observed test statistics within the group.  Differentially weighting
groups improves the power to detect signals in likely groupings.  The
advantage of the grouped-weighting concept, over fixed weights based
on prior information, is that it often leads to an increase in power
even if many of the groupings are not correlated with the signal.
Being data dependent, the procedure is remarkably robust to poor
choices in groupings. Power is typically improved if one (or more) of
the groups clusters multiple tests with signals, yet little power is
lost when the groupings are totally random.  If there is no apparent
signal in a group, relative to a group that appears to have several
tests with signals, the former group will be down-weighted relative to
the latter.  If no groups show apparent signals, then the weights will
be approximately equal.  The only restriction on the procedure is that
the number of groups be small, relative to the total number of tests
performed.

\vskip .05in


\newpage

Thorough testing for association between genetic variation and a
complex disease typically requires scanning large numbers of genetic
polymorphisms.  In a multiple testing situation, such as a whole
genome association scan, the null hypothesis is rejected for
any test that achieves a p-value less than a predetermined threshold.
To account for the greater risk of false positives, this threshold is
more stringent as the number of tests conducted increases. To bolster
power, recent statistical methods suggest up-weighting and
down-weighting of hypotheses, based on prior likelihood of association
with the phenotype (Genovese et al.  2006, Roeder et al. 2006).
Weighted procedures multiply the threshold by the weight $w$, for each
test, raising the threshold when $w>1$ and lowering it if $w<1$.  To
control the overall rate of false positives, a budget must be imposed
on the weighting scheme.  Large weights must be balanced with small
weights, so that the average weight is one.  These investigations
reveal that if the weights are informative, the procedure improves
power considerably, but, if the weights are uninformative, the loss in
power is usually small.  Surprisingly, aside from this budget
requirement, any set of non-negative weights is valid (Genovese et al.
2006).  While desirable in some respects, this flexibility makes it
difficult to select weights for a particular analysis.

The type of prior information readily available to investigators is
often non-specific.  For instance, SNPs might naturally be
grouped, based on features that make various candidates more
promising for this disease under investigation.  For a brain-disorder
phenotype we might cross-classify SNPs by categorical variables such
as those displayed in Table I.  The SNPs in $\calG_1$ seem most
promising, a priori, while those in $\calG_4$ seem least promising.
Those in $\calG_2$ and $\calG_3$ are more promising than those in
$\calG_4$, but somewhat ambiguous.  It is easy to imagine additional
variables that further partition the SNPs into various classes that
help to separate the more promising SNPs from the others.  While this
type of information lends itself to grouping SNPs, it does not lead
directly to weights for the groups.  Indeed it might not even be to
possible to choose a natural ordering of the groups. What is needed is
a way to use the data to determine the weights, once the
groups are formed.
\vskip .4in
\begin{center}
\begin{tabular}{ccc}\hline
&Functional & Non-Functional \\ \cline{2-3}
Brain expressed & $\calG_1$ & $\calG_2$\\
Non-Brain expressed & $\calG_3$ & $\calG_4$ \\ \hline
\end{tabular}
\end{center}
\noindent{Table I.}
\vskip .4in 

Until recently, methods for weighted multiple-testing
required that prior weights be developed independently of the data
under investigation (Genovese et al. 2006, Roeder et al. 2006).  In
this article we ask the following questions: if the weights are to be
applied to tests grouped by prior information, what choice of weights
will optimize the average power of the genetic association study?  How
can we estimate these weights from the data to achieve greater power
without affecting control of the family-wise error rate?  

\section*{Methods}

Consider
$m$ hypotheses corresponding to standardized test
statistics $T = (T_1, \ldots, T_m)$. The p-values associated with the tests are $(P_1, \ldots, P_m)$.
We assume $T_j$ is approximately
normally distributed with non-centrality parameter $\xi_j$, or
the tests are $\chi^2$ distributed with
non-centrality parameter $\xi_j^2$. When
using a Bonferroni correction for $m$ tests, the threshold for
rejection is achieved if the p-value $P_j \leq \alpha/m$.   
The
weighted Bonferroni procedure of Genovese, Roeder and Wasserman (2005)
is as follows.  Specify nonnegative weights $w=(w_1, \ldots, w_m)$ and
reject hypothesis $H_j$ if
\begin{equation}
\label{eq:reject}
j\in {\cal R} = \left\{ j:\ \frac{P_j}{w_j} \leq \frac{\alpha}{m}\right\}.
\end{equation}
As long as $m^{-1}\sum_j w_j =1$, this procedure controls family-wise
error rate at level $\alpha$.  For a test of $\xi_j = 0 \mbox{ vs. }
\xi_j \neq 0$, the power of a single weighted test is
$$
\pi(\xi_j,w_j)= \overline{\Phi}
\left(\overline{\Phi}^{-1}\left(\frac{\alpha w_j}{2m}\right)-\xi_j\right) +
\overline{\Phi}
\left(\overline{\Phi}^{-1}\left(\frac{\alpha w_j}{2m}\right)+\xi_j\right),
$$
where $\overline\Phi(t)$ is the upper tail probability of a
standard normal cumulative distribution function.  When the
alternative hypothesis is true, weighting increases the power when
$w_j >1$ and decreases the power when $w_j < 1$.  We call
$\pi(\xi_j,w_j)$ the {\em per-hypothesis power}.  For signals $(\xi_1,
\ldots, \xi_m)$ and weights $(w_1, \ldots, w_m)$ the {\em average
  power} is
\[
\overline{\pi}(\theta,w) = \frac{1}{m_1}\sum_{j=1}^m \pi(\xi_j,w_j).
\]

The optimal weight vector $w=(w_1, \ldots, w_m)$ that maximizes the
average power subject to $w_j \geq 0$ and $m^{-1}\sum_{j=1}^m w_j=1$
is (Wasserman and Roeder 2006)
\begin{equation}
\label{eq:wq1}
w(\xi_j) = 
\frac{m}{\alpha}\overline{\Phi}\left(\frac{|\xi_j|}{2} + 
\frac{c}{|\xi_j|}\right),
\end{equation}
where $c$ is the constant that satisfies the budget criterion on weights
\begin{equation}
\label{eq:wbest}
\frac{1}{m} \sum_{j=1}^m w(\xi_j)  =1.
\end{equation}

The optimal weights vary with the signal strength in a
non-monotonic manner (Figure \ref{fig::plotc}).  
For any particular sample, $c$ adjusts the weights to satisfy the
budget constraint on weights.  In so doing, it shifts the mode of the
weight function from left to right depending on the number of small,
versus large, signals observed.

The optimal weight function has an interesting effect
on the rejection threshold.  This choice of weights results in a
threshold for rejection that varies smoothly with the signal strength.
Figure \ref{fig::threshold} plots the rejection threshold
$-\log_{10}\,(\alpha w_j/m)$, calculated for the data displayed in Figure
\ref{fig::plotc}, as a function of the signal strength and
contrasts it with the rejection threshold of a Bonferroni corrected
test $-\log_{10}\,(\alpha/m$).  From Figures
\ref{fig::plotc}-\ref{fig::threshold} it is evident why an optimally weighted
test has greater power than a non-weighted test.  The weighted-threshold is
less stringent for signals in the midrange, and more stringent for
both large and small signals.  Consequently, if the signal is likely
to be very strong or very weak, the test is down-weighted (weight less
than one).  In practice, little power is lost by this
tradeoff. For small signals the chance of rejecting the hypothesis is
minimal with or without weights.  For large signals the p-value is
likely to cross the threshold regardless of the weight.  Larger
weights are focused in the midrange to help to reveal signals that are
marginal.

Clearly $\xi_j$ is not known, so it must be estimated to utilize this
weight function.  A natural choice is to build on the two stage
experimental design (Satagopan and Elston RC 2003; Wang et al. 2006)
and split the data into subsets, using one subset to estimate $\xi_i$,
and hence $w(\xi_i)$, and the second to conduct a weighted test of the
hypothesis (Rubin et al. 2006).  This approach would arise naturally
in an association test conducted in stages.  It does lead to a gain in
power relative to unweighted testing of stage 2 data; however, it is
not better than simply using the full data set without weights for the
analysis (Rubin et al. 2006; Wasserman and Roeder 2006). These results
are corroborated by Skol et al. (2005) in a related context.  They
showed that it is better to use stages 1 and 2 jointly, rather than
using stage 2 as an independent replication of stage 1.

To gain a strong advantage with data-based weights, prior information
is needed.  One option is to order the tests (Rubin et al. 2006), but
with a large number of tests this can be challenging.  Another option
is to group tests that are likely to have a signal, based on prior
knowledge, as follows:
\begin{enumerate}
\item Partition the tests into subsets $\calG_1,\ldots,\calG_K$, with the
  $k$'th group containing $r_k$ elements, ensuring that $r_k$ is at
  least 10-20.
\item Calculate the sample mean $Y_k$ and variance $S_k^2$ for the
  test statistics in each group.
\item Label the $i$'th test in group $k$, $T_{ik}$.  At
  best only a fraction of the elements in each group will have a signal, hence we
  assume that for $i=1,\ldots,r_k$ the distribution of the test
  statistics is approximated by a mixture model
\[ T_{ik} \sim (1-\pi_k) N(0,1) + \pi_k N(\xi_k,1)\]
or
\[ T_{ik} \sim (1-\pi_k) \chi_1^2(0) + \pi_k \chi_1^2(\xi_k^2)\]
where  $\xi_{k}$ is the signal size for those tests with a signal in the $k$'th group. (This is an approximation because the signal is likely to vary across tests.)
\item Estimate $(\pi_k,\xi_k)$ using the method of moments 
estimator.  For the normal model this is
\[\hat \pi_k= Y_k^2/(Y_k^2 + S_k^2 - 1),\qquad
\hat \xi_k = Y_k/\pi_k,\] provided $\hat \pi_k > 1/r_k$; otherwise $\hat \xi_k=0$.
\\

For the $\chi^2$ model $\hat \xi^2$ is a root of the quadratic
equation $x^2 -bx +1=0$ where $b = (S_k^2-1)/(Y_k-1) + Y_k -5$.
If both roots are negative, $\hat \xi_k^2=0$; otherwise,
$\hat \pi_k = (Y_k-1)/\hat \xi_k^2$.  

\item For each of the $k$ groups, construct weights $w(\hat \xi_k)$.
  Then, to account for excessive variability in the weights, induced
  by variability in $\hat \xi_k$, smooth the weights by taking a
\[\hat w_k = 0.95 \,w(\hat \xi_k) + 0.05\, K^{-1}\sum_k w(\hat \xi_k).\] 
Renorm weights if necessary to ensure the weights sum to
$m$.  Each test in group $k$ receives the
weight $\hat w_k$. 
\end{enumerate}

This weighting scheme relies on data-based estimators of the optimal
weights, but with a partition of the data sufficiently crude to
preserve the control of family-wise error rate.  The approach is an
example of the ``sieve principle''.  More formally this result is
stated in the following Theorem.   \vskip
.25in
\noindent {\bf Theorem.}
Let $b_m = \frac 1m\sum_k \sqrt{r_k}$.  If $\sum_{j=1}^m \hat w_j =m$,
then ${\cal R}$ (\ref{eq:reject}) controls family-wise error at level
$\alpha + O(b_m)$. Proof is in the Appendix.  \vskip .25in
\noindent This result establishes control of family-wise error 
at level $\alpha$,
asymptotically, provided
\[b_m = \frac {\sum_k \sqrt{r_k}}
{\sum_k {r_k}}\to 0,\qquad \mbox{ as } m\to\infty.\] The inflation
term in the error rate is near zero under a number of circumstances.
Loosely speaking, the requirement is that each group contains a
sufficient number of elements to permit valid estimation of $\{\hat
\xi_k\}$.  For instance, if each group has the same number of
elements $r_k = r$, then $b_m=1/\sqrt{r}$, which goes to zero, provided
the number of groups grows more slowly than the number of tests
performed.  Likewise, $b_m\to 0$ if ${\max\{\sqrt{r_k}\}}/
\min\{r_k\}\to 0$.

Figure \ref{fig::meanvar} illustrates how $w(\hat \xi_k)$ varies with
$\hat \xi_k$ and the sample variances (weight is proportional to the
diameter of the circle).  Notice that weight increases as a function
of the signal until it becomes fairly large and then declines.

\section*{Results}

To simulate a large scale study of association, we generate test
statistics from $m=10,000$ tests with $m_1=50$ and $100$ tests having
a signal ($\xi_i > 0$) and $m_0=m-m_1$ following the null hypothesis.
These choices were made to simulate the second stage of a two-stage
genome-wide association study, with about 1/3-1\% of the initial SNPs
tested at stage 2. In the proximity of a causal SNP, clusters of tests
tend to exhibit a signal. We simulate the data as if 5-10 additional
SNPs were in the proximity of each causal SNP.  Thus, if 10-20 actual
causal variants are present in the genome, approximately 50 to 100
tests might be associated with the phenotype at varying levels of
intensity.

The simulated signal strengths vary over 5 levels
$(\xi_1,\ldots,\xi_5) = \xi_0\times (1,1.5,2,2.5,3)$ 
with $m_1/5$ realizations of each of the
5 levels of signals. The $m$ simulated tests are grouped into categories 
$\calG_1,\ldots,\calG_K$ with the groupings formed to convey various levels of
informativeness. Let $\xi_{ik}$ be the signal of the $i$'th element in
group $k$, $\bar \xi_{.k}$ be the mean in group $k$, and $\bar
\xi_{..}$ be the mean of the whole set, respectively.  The information
in a prior grouping is summarized by the $R^2$
\[ R^2 = 1-\frac{\sum_k \sum_i (\xi_{ik}- \bar \xi_{.k})^2}   
{\sum_k \sum_i (\xi_{ik}- \bar \xi_{..})^2}.\] 


The $10,000$ tests are grouped into 10 categories.  We start the
process by dividing the $m_0$ tests that do not have a signal randomly
into 5 equal sized groupings, $\calG_1,\ldots,\calG_5$. Now $m_1$
tests remain to constitute the remaining 5 categories,
$\calG_6,\ldots,\calG_{10}$.  We create the ideal partition of these
tests by placing all tests with a common value of $\xi_j$ in the same
category. 
Next, to create more realistic groupings, we move some tests from
categories 1-5 into 6-10 and vice versa.  Specifically, we move a
fraction $p_0$ of the $m_0$ null tests to categories 6-10, and
distribute them evenly.  Likewise we move a fraction $p_1$ of the $m_1$
tests with $\xi>0$ to categories 1-5, and distribute them evenly.  By
varying $(p_0, p_1)$ we obtain various levels of informativeness of
the groupings, reflecting priors of various value.

To see the effect of including null loci in the same grouping as the
SNPs with true effects, we fix $(\xi_0=2,p_1=0,m_1=100)$ and vary
$p_0$.  Setting $p_0=0.5$ (0.1) increases the elements of groups 6-10
to 1,010 (218), but only 20 are true alternatives.  For $p_0=$ 0.01,
0.1, 0.25, and 0.5 we find a difference in power (weighted minus the
unweighted procedure) of 14, 5, 0, and -3 percent, respectively.  So,
for $p_0> 0.25$ there is a loss in power, but it is relatively small.

Next we explore the effect of failing to place the true effects in the
more promising categories (6-10). To do so, we fix
$(\xi_0=2,p_0=.1,m_1=100)$ and vary $p_1$.  For $p_1=$ 0.05.  0.1,
0.5, and 0.9, we find a difference in power of 7, 3, 2, -5 and -2
percent, respectively.  Even when 90\% of the true alternatives are
grouped with large numbers of nulls in groups 1-5, the loss in power
is relatively small.  Another interesting feature is that a 50\% swap
leads to a greater loss in power than a 90\% swap.  The latter occurs
because weights are approximately constant across groups when the
alternatives are scattered nearly at random.  When half of the
alternatives are in the promising groups, these categories are
up-weighted at the expense of the other categories.  This balance can
lead to a net loss in power, relative to the unweighted test.

Figure \ref{fig::powerR2} displays the difference in power as a
function of $R^2$.  The proportion of null tests in cells 1-5, and
alternative tests in cells 6-10 varies: $p_0\in [0.01-0.5]$ and
$p_1\in [0.01-0.95]$.  From these simulations we see that, provided
$p_0<0.5$ and $p_1>0.1$, the weighted method is generally more
powerful than the unweighted method (plot symbol ``o'').  Two
exception occur; both have $R^2$ less than 2\% of the variability in
signal.  For $R^2$ near 0 the loss in power from poorly selected
groupings is modest. Deviations in $p_1$ from ideal have a greater
impact than deviations of $p_0$ (plot symbol ``$\star$'' vs. ``$+$'').
This asymmetry is expected because groups (1-5) contain many more
elements than groups 6-10.  Consequently signals can be swamped by
nulls in these groupings.
Finally we tried mixing the various levels of true
alternatives $\xi_0\times (1,1.5,2,2.5,3)$ among groups 6-10 and
found that this had a negligible effect on the power (results not
shown).

\section*{Discussion}

Whole genome analysis has generated much discussion about power, the
effect of multiple testing on power, and various multistage
experimental designs (e.g., Wang et al. 2006). We investigate the
performance of a weighting scheme that allows for the input of weak
prior information, in the form of groupings of tests, to improve power
in large scale investigations of association.  The method can be
applied at any stage of an experiment. The beauty of the
grouped-weighting concept is that it is likely to lead to an increase
in power, provided multiple tests with signals are clustered together
in one (or more) of the groups. Little power is lost when many groups
contain no true signal. This remarkable robustness is achieved because
the procedure uses the observed test statistics in the grouping to
determine the weight.  If there is no apparent signal, the group will
be down-weighted.  The only restriction on the procedure is that the
number of groups be small, relative to the total number of tests
performed.

Using groupings and weights to interpret the many tests conducted in a
large scale association study has potential, regardless of power lost
when weights are poorly chosen. 
Typically some SNPs are favored 
due to knowledge gleaned from the literature and prior
investigations.  When seemingly random SNPs produce smaller p-values
than the favored candidates, one is baffled about how to handle the
situation.  Moreover, it often happens that promising candidate SNPs
do produce small p-values, but these p-values might not be small enough to
cross the significance threshold when a Bonferroni correction is applied.
After the huge investment of a whole genome scan it would
be foolhardy not to pursue both (i) SNPs that produce tiny p-values
and (ii) SNPs that produce respectable p-values that would have been
significant had a formal weighting scheme been utilized to incorporate
prior information.  We suggest using the weighting method of analysis
described here as a way to formalize the incorporation of prior
information.  

Weights can be incorporated into various multiple testing procedures,
including false discovery methods.  This paper considers
controlling family-wise error rate, but similar results hold for false
discovery control (Benjamini and Hochberg 1995) and will be pursued
elsewhere.

\section*{Appendix}

\noindent
{\sf Proof of Theorem 1.}  Let ${\cal H}_0$ denote
the set of indices for which $\xi_j=0$.  With fixed weights, the
family-wise error is
\begin{eqnarray*}
\P(({\cal R}\cap {\cal H}_0) > 0)
  &  =   & \P\left( P_j \le \frac{\alpha w_j }{m}\ \ {\rm for\ some\ }j \in {\cal H}_0\right)\\
  & \leq & \sum_{j\in {\cal H}_0} \P\left( P_j \le \frac{\alpha w_j }{m}\right)=
  \frac{\alpha}{m} \sum_{j\in {\cal H}_0} w_j \leq \alpha \overline{w} = \alpha .\ \ \ 
\end{eqnarray*}
The estimated signal in the group occupied by the $j$'th test, $\hat
\xi_k$ is estimated from a sample of $r_k$ test statistics,
consequently $\hat \xi_k = \xi_k + O\left(r_k^{-1/2}\right)$.  Thus
with random weights
\begin{eqnarray*}
\P(({\cal R}\cap {\cal H}_0) > 0)
  &\leq & \sum_{j\in {\cal H}_0} \P\left( P_j \le \frac{\alpha w_j(\hat \xi_k) }{m}\right)\\
  & \approx   & \frac{\alpha}{m} \sum_{j\in {\cal H}_0}\left\{ w_j(\xi_k)+\left(w_j(\hat\xi_k)-w_j(\xi_k)\right) \right\}\\  
  & \leq & \alpha(1+O(b_m)). 
\end{eqnarray*}

\vspace{1cm}

\vspace{1cm}

\section*{References}

\begingroup
\parindent=0em
\baselineskip=17pt
\parskip=6pt
\everypar={\hangindent=1em\hangafter=1\noindent}

Benjamini Y, Hochberg Y (1995).
Controlling the false discovery rate: a practical and powerful
approach to multiple testing. 
J Roy Stat Soc B 57:289-300

Genovese, CR, Roeder, K, Wasserman, L (2006).
False Discovery Control with p-Value Weighting.
Biometrika 93:509-524.

Roeder, Bacanu, Wasserman and Devlin (2005).
Using Linkage Genome Scans to
Improve Power of Association in Genome Scans.
{\em The American Journal of Human Genetics.}
{\bf 78}.

Rubin, D, van der Laan, M. and Dudoit, S. (2006). 
Multiple testing procedures which are optimal at a simple alternative. 
Collection of Biostatistics Research Archive, 
http://www.bepress.com/ucbbiostat/paper171/

Satagopan JM, Elston RC (2003)       
Optimal two-stage genotyping in population-based association studies.
Genet Epidemiol  25:149-57

Skol AD, Scott LJ, Abecasis GR, Boehnke M (2006) 
Joint analysis is more efficient than replication-based 
analysis for two-stage genome-wide association studies.
Nat Genet. 38:209-213.

Wang H, Thomas DC, Pe'er I, Stram DO.  
Optimal two-stage genotyping designs for genome-wide association scans.
Genet Epidemiol. 2006 May;30(4):356-68. 

Wasserman L, Roeder K, (2006) Weighted Hypothesis Testing.\\
http://arxiv.org/abs/math.ST/0604172

\begin{figure}
\begin{center}
\includegraphics[angle=-90,width=4in]{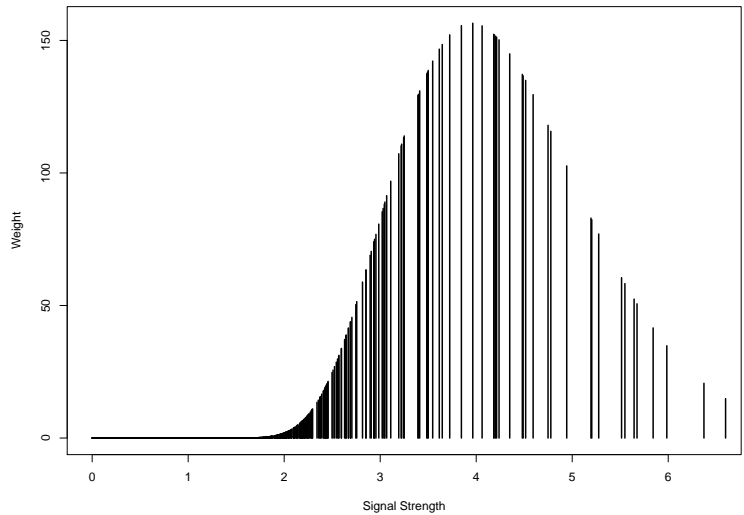}
\end{center}
\caption{Distribution of optimal weights for $m=100,000$ simulated 
tests (a random selection of
  5,000 are displayed).  The signal strength is the non-centrality
  parameter for a standard normal test statistic; if the test
  statistic is $\chi^2$ distributed, the signal strength is the square
  root of the non-centrality parameter.  }
\label{fig::plotc}
\end{figure}

\begin{figure}
\begin{center}
\includegraphics[angle=-90,width=4in]{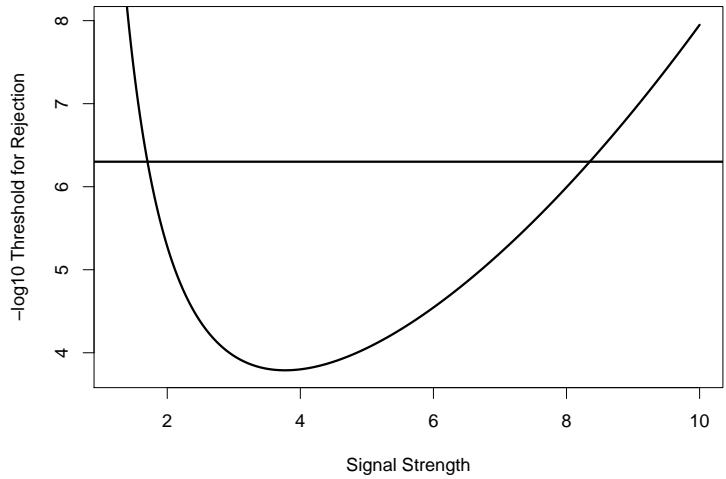}
\end{center}
\caption{Threshold for rejecting P-values versus signal strength.  The 
  $log_{10}$ p-value is rejected if it is larger than the threshold.
  For this illustration $m=100,000$ and $\alpha=0.05$.  The unweighted
  Bonferroni has a constant threshold value (horizontal line).  The
  weighted threshold varies as a function of the weight (curved line).
  The optimal weight is calculated as a function of the (estimated)
  signal strength.}
\label{fig::threshold}
\end{figure}

\begin{figure}
\begin{center}
  \includegraphics[angle=-90,width=4in]{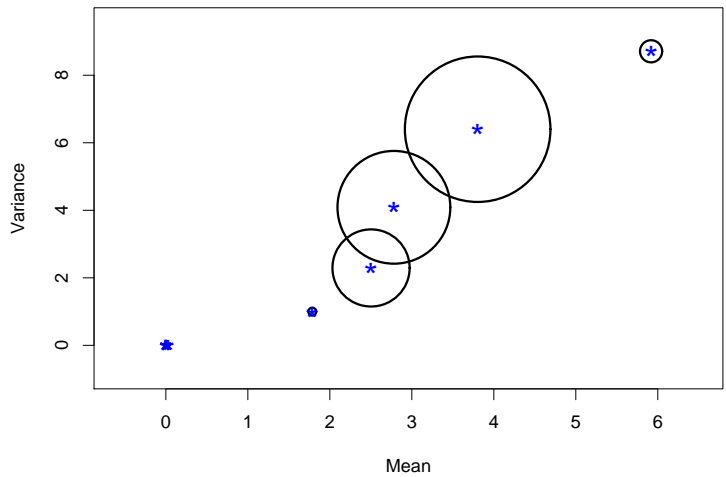} \end{center}
\caption{Weight as a function of $\hat\xi_k$ and variance.  The diameter of the circle indicates relative weight.}  
\label{fig::meanvar}
\end{figure}

\begin{figure}
\begin{center}
  \includegraphics[angle=-90,width=4in]{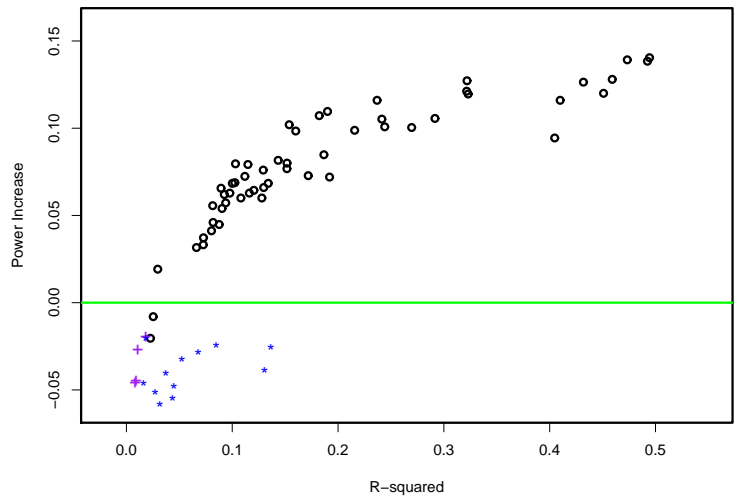} 
\end{center}
\caption{Net power different between weighted Bonferroni and
  unweighted, as a function of $R^2$. The worst cases are
$p_0 = 0.5$ (plot symbol +) and $p_1>0.1$ (plot symbol *). The remaining models have plot symbol o.}
\label{fig::powerR2}
\end{figure}

\end{document}